\documentclass[11pt,letterpaper,reqno]{amsart}
\usepackage{amsmath,amsthm,amsfonts,amssymb,mathtools,algpseudocode}
\usepackage{comment,bm,color,mathdots,subfig,bookmark}

\hypersetup{pdfstartview={FitH}}

\newtheorem{Theorem}{{\bf Theorem}}[section]

\newcommand{\diag}{\textnormal{diag}}

\newcommand{\im}{\textnormal{Im}}
\newcommand{\re}{\textnormal{Re}}
\newcommand{\N}{\mathcal{N}}
\newcommand{\U}{\mathcal{U}}
\newcommand{\Ham}{\mathcal{H}}
\newcommand{\sHam}{\mathcal{W}}
\newcommand{\perH}{\mathcal{M}}
\newcommand{\sperH}{\mathcal{K}}
\newcommand{\Symp}{\mathcal{S}p}
\newcommand{\Perp}{\mathcal{P}p}
\numberwithin{equation}{section}

\begin{document}

\title[On Normal Structured Matrices]{On normal and structured matrices under unitary structure-preserving transformations}

%\author{One author\thanks{some info} \and Another author\thanks{more info}}
\author{Erna Begovi\'{c}~Kova\v{c}}\thanks{\textsc{Erna Begovi\'{c}~Kova\v{c},
Faculty of Chemical Engineering and Technology, University of Zagreb, Maruli\'{c}ev trg 19, 10000 Zagreb, Croatia}.
\texttt{ebegovic@fkit.hr}}
\author{Heike Fa{\ss}bender}\thanks{\textsc{Heike Fa{\ss}bender,
Institute for Numerical Analysis, TU Braunschweig, Universit\"atsplatz 2, 38106 Braunschweig, Germany}.
\texttt{h.fassbender@tu-braunschweig.de}}
\author{Philip Saltenberger}\thanks{\textsc{Philip Saltenberger,
Institute for Numerical Analysis, TU Braunschweig, Universit\"atsplatz 2, 38106 Braunschweig, Germany}.
\texttt{philip.saltenberger@tu-braunschweig.de}}

%\thanks{}
\date{\today}


\begin{abstract}
Structured canonical forms under unitary and suitable structure-preserving similarity transformations for 
normal and (skew-)Hamiltonian as well as normal and per(skew)-Hermitian matrices are proposed.
Moreover, an algorithm for computing those canonical forms is sketched.

\end{abstract}

\maketitle

\section{Introduction}
It is well known that  the set $\N_n$  of  normal matrices 
\[\N_n=\{X \in \mathbb{C}^{n \times n}\ : \ XX^H=X^HX\}\]
can be described as the set of matrices that are unitarily
diagonalizable, that is,
\[ \N_n = \{ X \in \mathbb{C}^{n \times n}\ : \text{ there exists } U \in \U_n \text{ such that } U^HXU \text{ is diagonal}\}
\]
for
\[
\U_n = \{ U \in \mathbb{C}^{n \times n}\ :  U^HU = UU^H = I \}.
\]

Here we will consider normal matrices $A \in \N_{2n}$ which belong to the following sets $\mathcal{S}$  of structured matrices:
\begin{align*}
\Ham &= \{ H \in \mathbb{C}^{2n \times 2n} \ : (JH)^H = JH\}  &\textnormal{ Hamiltonian matrices,}\\
\sHam &= \{ W \in \mathbb{C}^{ 2n \times 2n} \ : (JW)^H = -JW\}  &\textnormal{ skew-Hamiltonian matrices,}\\
\perH &= \{ M \in \mathbb{C}^{ 2n \times 2n} \ : (FM)^H = FM\}  &\textnormal{ per-Hermitian matrices,}\\
\sperH &= \{ K \in \mathbb{C}^{ 2n \times 2n} \ : (FK)^H = -FK\}  &\textnormal{ perskew-Hermitian matrices,}
\end{align*}
where
$$J = J_{2n}=\left[
      \begin{array}{cc}
        0 & I_n \\
        -I_n & 0 \\
      \end{array}
    \right] \in \mathbb{R}^{2n \times 2n} , \quad
F_{n}=\left[
\begin{array}{ccc}
&&1\\
&\iddots&\\
1&&
  \end{array}
\right]\in \mathbb{R}^{n \times n}\quad \textnormal{and} \quad F = F_{2n}.$$
Certainly, any normal and structured matrix $A \in \N_{2n} \cap {\mathcal S},  \mathcal{S} \in \{\Ham, \sHam, \perH, \sperH\}$ is unitarily diagonalizable by some $U \in \U_{2n},$ but  in general the matrix $U^HAU$ is not in $\mathcal{S}.$ 
 Typically, $U$ will have to have an additional property in order to force $U^HAU \in \mathcal{S}.$ 
It is well-known that
 transformations that preserve the structure of the sets $\Ham$ and $\sHam$ are symplectic transformations, while  transformations that preserve the structure of
 the other two sets are perplectic transformations, where
\begin{align*}
\Symp &= \{ S \in \mathbb{C}^{2n \times 2n} \ : S^HJS=J\}  &\textnormal{symplectic matrices,}\\
\Perp &= \{ P \in \mathbb{C}^{ 2n \times 2n} \ : P^HFP = F\}  &\textnormal{perplectic matrices}.
\end{align*}
That is, $Z^HAZ \in \Ham$ (resp. $\sHam$) for $A \in \Ham$ (resp. $\sHam$) and $Z \in \Symp,$
and $Z^HAZ \in \perH$ (resp. $\sperH$) for $A \in \perH$ (resp. $\sperH$) and $Z \in \Perp.$
Numerous properties of the sets $\Ham, \sHam, \perH, \sperH, \Symp, \Perp$ (and their interplay)
have been studied in the literature, see, e.g., \cite{MMT03,Trench04} and the references therein.

Here, we are interested in unitary symplectic and unitary perplectic transformations,
\begin{align*}
\U_{2n} \cap \Symp & = \{ S \in \mathbb{C}^{2n \times 2n} \ : S = \begin{bmatrix} S_1 & S_2\\ -S_2 & S_1\end{bmatrix},
S^HS = I_{2n}, S_1, S_2 \in \mathbb{C}^{n \times n} \},\\
\U_{2n} \cap \Perp &= \{ P \in \mathbb{C}^{2n \times 2n} \ : P^HFP=F, P^HP=I_{2n} \}.
\end{align*}
Even so any matrix $A \in \N_{2n}\cap \mathcal{S}, \mathcal{S} \in \{\Ham, \sHam, \perH, \sperH\}$ is unitarily diagonalizable, there need not be a matrix $Z \in \U_{2n} \cap \Symp ~(\text{resp., } \U_{2n} \cap \Perp)$  such that $Z^HAZ$ is diagonal.
For each of the four sets $\Ham, \sHam, \perH, \sperH$ of matrices, we will 
give the canonical forms for matrices $A \in \N_{2n} \cap {\mathcal S},$ $\mathcal{S} \in \{\Ham, \sHam, \perH, \sperH\}$ under unitary structure-preserving transformations.
For the (skew-)Hamiltonian case, the structured canonical form for normal matrices is deduced from the Hamiltonian
Schur form \cite{LinMX99,PaiVD86}. As no structured canonical form for normal per(skew-)Hermitian matrices
is known, their structured canonical form is derived directly.

%%%%%%%%%%%%%%%% 
Goldstine and Horwitz describe in \cite{GolH59} an approach for diagonalizing a normal matrix  
$A\in \mathbb{C}^{n \times n}.$  It is based on the fact that any matrix $A\in \mathbb{C}^{n \times n}$
can be written as the sum of its Hermitian part $B= \frac{1}{2}(A + A^H)$ and its skew-Hermitian
part $C= \frac{1}{2}(A - A^H).$ In a first step, the Hermitian matrix $B$ is diagonalized (e.g., 
the classical Jacobi method \cite{Jac1846,GolVL12} adapted for Hermitian matrices \cite{Hen58}).
Applying the necessary transformations not just to $B$, but to $C$ as well will not inevitably transform
 $C$ to diagonal form as well. But the nonzero off-diagonal elements of the resulting matrix can easily
be eliminated by Givens transformation without wrecking the diagonal form of $B.$
Based on these ideas we will give an algorithm for computing the structured canonical form
of normal (skew-)Hamiltonian matrices. That algorithm will make use of the 
 Jacobi-type algorithm for computing the Hamiltonian Schur form of a Hamiltonian matrix
with no purely imaginary eigenvalues given in \cite{BunF97}.
It can be applied to diagonalize a Hermitian Hamiltonian matrix,
but it can not be applied to a normal Hamiltonian matrix with purely imaginary eigenvalues. 

In Section~\ref{sec:structured_ham} the structured canonical form for normal Hamiltonian and normal skew-Hamil\-tonian matrices is presented,
while in Section~\ref{sec:structured_perp} the per-Hermitian and perskew-Hermitian case is dealt with.
Section~\ref{sec:Jac_Ham} adapts the approach from \cite{GolH59} for diagonalizing a normal matrix  
 to computing the structured canonical form of normal (skew-)Hamiltonian matrices. Some of its numerical properties
will be discussed.
In a similar fashion 
 algorithms for computing the canonical form of normal per-Hermitian and normal perskew-Hermitian matrices 
can be derived. The details are omitted here. 

%%%%%%%%%%%%%%%%%%%%%%%%%%%%%%%%%%%%%%%%%%%%%
\section{Structured canonical form for normal (skew-)Hamiltonian matrices under unitary symplectic similarity transformation}\label{sec:structured_ham}

The (complex) Hamiltonian matrices\footnote{These matrices have been called  $J$-Hermitian matrices
in \cite{MMT03}, in order to distinguish them from the matrices $H \in \mathbb{C}^{2n \times 2n}$ which satisfy $(JH)^T=JH$. The latter ones are called $J$-symmetric in \cite{MMT03}. As there is no ambiguity
here, we will simply use the term Hamiltonian.}
form a Lie algebra and the (complex) skew-Hamiltonian matrices\footnote{These matrices have been called  $J$-skew-Hermitian matrices
in \cite{MMT03}, in order to distinguish them from the matrices $W \in \mathbb{C}^{2n \times 2n}$ which satisfy $(JW)^T=-JW$. The latter ones are called $J$-skew-symmetric in \cite{MMT03}. As there is no ambiguity
here, we will simply use the term skew-Hamiltonian.} form a Jordan algebra associated with the skew-Hermitian sesquilinear form $x^HJy.$ Both classes of matrices are well studied, see, e.g., \cite{MMT03} and the references therein. The symplectic matrices\footnote{These matrices have been called conjugate symplectic in \cite{MMT03}, in order to distinguish them from the matrices $S \in \mathbb{C}^{2n \times 2n}$ which satisfy $ S^TJS=J$. The latter ones are called complex symplectic in \cite{MMT03}. As there is no ambiguity
here, we will simply use the term symplectic.} form the automorphism group associated with the skew-Hermitian sesquilinear form $x^HJy.$

\subsection{Canonical form for normal Hamiltonian matrices }\label{sec:subsec1}
A matrix  $H\in\mathbb{C}^{2n\times2n}$ is (complex) Hamiltonian if $(JH)^H=JH,$
or equivalently, $H^H=JHJ$.
If we write $H$ as a $2\times2$ block matrix of $n \times n$ blocks, it is easy to verify that
\begin{equation}\label{eq:Hamsetup}
H=\left[
      \begin{array}{cc}
        H_{11} & H_{12} \\
        H_{21} & -H_{11}^H \\
      \end{array}
    \right], \qquad \textnormal{where} \ \ H_{12}^H=H_{12}, \ H_{21}^H=H_{21}.
\end{equation}
The eigenvalues of Hamiltonian matrices $H$ come in pairs $(\lambda, -\bar{\lambda})$
with $\lambda$ and $-\bar{\lambda}$ having the same multiplicity: if $\lambda \in \sigma(H),$ then
 $\bar{\lambda} \in \sigma(H^H) = \sigma(JHJ) = \sigma(-J^{-1}HJ) = -\sigma(H),$ where $\sigma(X)$ denotes the spectrum of the square matrix $X.$

Obviously, a diagonal Hamiltonian matrix $\hat{H}$ is of the form
\begin{equation}\label{eq:diagHam}
\hat{H}=\left[
       \begin{array}{cc}
         \Lambda & 0 \\
         0 & -\Lambda^H \\
       \end{array}
     \right]
\end{equation}
with a diagonal matrix $\Lambda = \diag(\lambda_1, \ldots, \lambda_n).$

A normal Hamiltonian matrix $H$ can be unitarily diagonalized;
\[ U^HHU = \diag(\mu_1, \ldots, \mu_{2n}) =:D\]
with a unitary matrix $U.$ Clearly, $\mu_i \in \sigma(H) = \{ \lambda_1, \ldots, \lambda_n, -\bar{\lambda}_1, \ldots, -\bar{\lambda}_n\}.$ The eigenvalues of $H$ can appear in any order on the diagonal $D$.
In particular, $U$ can be chosen such
that $D=U^HHU$ is a diagonal Hamiltonian matrix as in \eqref{eq:diagHam}. But, this $U$ will not be symplectic in general. 

When we restrict our transformations to unitary and symplectic ones, then we might not be able to diagonalize a normal
Hamiltonian matrix. Consider the Hamiltonian matrix $J_2 = \left[\begin{smallmatrix} 0 & 1\\ -1 & 0\end{smallmatrix}\right].$
It is normal, but any symplectic transformation does not change it at all. Thus we will not be able to find a unitary and
symplectic transformation which diagonalizes $J_2.$

The most condensed form for normal and Hamiltonian matrices under unitary and symplectic transformation
which can be achieved is derived next.
We will start with the Hamiltonian Schur form presented in \cite{LinMX99} (see also \cite{PvL81}).
\begin{Theorem}[Hamiltonian Schur form] \label{thm:Ham}
For any Hamiltonian matrix $H \in \mathbb{C}^{2n \times 2n}$ there exists a unitary and symplectic
matrix $U$ such that
\[ U^HHU = \begin{bmatrix}
T_{11} & T_{12} & X_{11} &X_{12}\\
0 & T_{22} & X_{21} & X_{22}\\
0 & 0 & -T_{11}^H & 0 \\
0 & Y_{22} & -T_{12}^H & -T_{22}^H
\end{bmatrix} \in \Ham
\]
where $T_{11}\in \mathbb{C}^{n_1 \times n_1}$ is upper triangular and $\left[\begin{smallmatrix} T_{22}& X_{22}\\Y_{22} & -T_{22}^H\end{smallmatrix}\right]\in \mathbb{C}^{2n_2 \times 2n_2}$ is a Hamiltonian matrix with purely imaginary eigenvalues; $n_1 + n_2 = n.$
\end{Theorem}
If the Hamiltonian matrix $H$ is normal, that is, if $HH^H=H^HH$, then it follows from Theorem \ref{thm:Ham} that
\[ U^HHU = \begin{bmatrix}
D_{1} & 0 & 0 & 0\\
0 & T_{22} & 0 & X_{22}\\
0 & 0 & -D_{1}^H & 0 \\
0 & Y_{22} & 0 & -T_{22}^H
\end{bmatrix}
\]
where $D_{1}$ is diagonal and $\hat{H}=\left[\begin{smallmatrix} T_{22}& X_{22}\\Y_{22} & -T_{22}^H\end{smallmatrix}\right] $ is a normal Hamiltonian matrix with purely imaginary eigenvalues. In particular, it holds $X_{22}^H = X_{22}.$
Moreover, a normal matrix with purely imaginary eigenvalues  is skew-Hermitian. This implies
$\hat{H} = -\hat{H}^H.$ From this we have $T_{22} = -T_{22}^H$ and $Y_{22} = -X_{22}^H.$ Thus,
\[
\hat{H} = \begin{bmatrix}
T_{22}& X_{22}\\-X_{22} & T_{22}
\end{bmatrix}, \qquad X_{22} = X_{22}^H, \quad T_{22} = -T_{22}^H.
\]
Now let $Q= \frac{1}{\sqrt2} \left[\begin{smallmatrix} I & \imath I\\ \imath I  &I\end{smallmatrix}\right].$
Q is unitary and block-diagonalizes $\hat{H},$
\[
M = Q^H \hat{H} Q=
\begin{bmatrix}
 T_{22}+\imath X_{22}  &0 \\
 0 & T_{22}-\imath X_{22}
\end{bmatrix}.
\]
$M$ is skew-Hermitian, $M^H = -M.$
Thus, the block matrices $T_{22}+\imath X_{22} $ and $T_{22}-\imath X_{22} $ are skew-Hermitian
and can be diagonalized by unitary matrices
$V_1$ and $V_2,$ respectively. Hence,  $V^HMV$ is diagonal with the unitary matrix $V = \diag(V_1,V_2).$
Finally, transforming $V$ by $Q,$ that is, $\hat{S} = QVQ^H$ yields a symplectic and unitary matrix $\hat{S}$
which diagonalizes the four blocks of $\hat{H},$
\begin{equation}\label{eq_4d}
\hat{S}^H\hat{H} \hat{S} = \left[\begin{smallmatrix} D_2& D_3\\-D_3 & D_2\end{smallmatrix}\right]
\end{equation}
with $D_2 = -D_2^H,  D_3=D_3^H.$  This implies that  $D_2$ is a diagonal matrix with purely imaginary diagonal elements, while
 $D_3$ is diagonal with real diagonal elements.
Now we partition $\hat{S} = \left[\begin{smallmatrix} S_{11}& S_{12}\\S_{21} & S_{22}\end{smallmatrix}\right]$ into four square blocks conformal to \eqref{eq_4d}
and embed $\hat{S}$ into a $2n \times 2n$ identity matrix,
\[
S = \left[ \begin{array}{cc|cc}
I_{n_1} &  & & \\
 & S_{11} & & S_{12}\\ \hline
&& I_{n_1} &\\
& S_{21} &&S_{22}
\end{array}\right],
\]
such that the four blocks are of size $n \times n.$ Then $US$ is a unitary and symplectic matrix which diagonalizes all
four $n \times n$ blocks of $H.$
In summary, we have
\begin{Theorem}[Hamiltonian Schur form for normal Hamiltonian matrices]\label{theo:HamSchurNormal}
 For any normal Hamiltonian matrix $H \in \mathbb{C}^{2n \times 2n}$ there exists a unitary and symplectic
matrix $Z$ such that
\begin{equation}\label{eq:HamCanonical1}
 Z^HHZ = \begin{bmatrix}
D_1 & 0 & 0 &0\\
0 & D_2 & 0 & D_3\\
0 & 0 & -D_1^H & 0 \\
0 & -D_3 & 0 & D_2
\end{bmatrix} \in \Ham
\end{equation}
where $D_j, j = 1, 2, 3$ are diagonal matrices, $D_1 \in \mathbb{C}^{n_1 \times n_1}, D_2 \in \imath \mathbb{R}^{n_2 \times n_2}, D_3 \in \mathbb{R}^{n_2 \times n_2}, n_1+n_2=n,$ and
$\left[\begin{smallmatrix} D_2& D_3\\-D_3 & D_2\end{smallmatrix}\right]$ is a Hamiltonian and skew-Hermitian matrix with purely imaginary eigenvalues.
\end{Theorem}
%Hence, the canonical form for normal and Hamiltonian matrices $H \in \N_{2n} \cap \Ham$ under a unitary and symplectic transformation $Z \in \U_{2n} \cap \Symp$ can be depicted as
%\begin{equation}\label{eq:HamCanonical1}
%Z^HHZ = \begin{bmatrix}
%\Diag & \Diag \\ \Diag & \Diag
%\end{bmatrix} =
%\begin{bmatrix}
%\Lambda_1 & \Lambda_2 \\ -\Lambda_2 & -\Lambda_1^H
%\end{bmatrix} =: \Lambda
%\end{equation}
%with diagonal matrices $\Lambda_1 \in \mathbb{C}^{n \times n}$ and $\Lambda_2 \in \mathbb{R}^{n \times n}.$
%Moreover, if the $i$th diagonal element of $\Lambda_2$ is nonzero, then the $i$th diagonal element of
%$\Lambda_1$ is purely imaginary (or zero).

%%%%%%%%%%%%%%%%%%%%%%%%%%

\subsection{Canonical form for normal skew-Hamiltonian matrices}
A matrix $W\in\mathbb{C}^{2n\times2n}$ is (complex) skew-Hamiltonian if
$(JW)^H=-JW,$ or equivalently, $W^H=-JWJ$.
%Written as a $2\times2$ block matrix of $n \times n$ blocks it takes the form
%$$W=\left[
%      \begin{array}{cc}
%        W_{11} & W_{12} \\
%        W_{21} & W_{11}^H \\
%      \end{array}
%    \right], \qquad W_{12}^H=-W_{12}, \ W_{21}^H=-W_{21}.$$
%The eigenvalues of skew-Hamiltonian matrices $W$ come in pairs $(\lambda, \bar{\lambda})$
%with $\lambda$ and $\bar{\lambda}$ having the same multiplicity: if $\lambda \in \sigma(W),$ then
% $\bar{\lambda} \in \sigma(W^H) = \sigma(-JWJ) = \sigma(J^{-1}WJ) = \sigma(W).$
%Obviously, a diagonal skew-Hamiltonian matrix $\hat{W}$ is of the form
%\begin{equation}\label{eq:diagskewHam}
%\hat{W}=\left[
%       \begin{array}{cc}
%         \Lambda & 0 \\
%         0 & \Lambda^H \\
%       \end{array}
%     \right]
%\end{equation}
%with a diagonal matrix $\Lambda = \diag(\lambda_1, \ldots, \lambda_n).$
%
%A normal skew-Hamiltonian matrix $W$ can be unitarily diagonalized. In particular, the unitary transformation
%matrix $U$ can be chosen such that $U^HWU$ is a diagonal skew-Hamiltonian matrix as in \eqref{eq:diagskewHam}.
%But, as in the Hamiltonian case, $U$ will not be symplectic in general.

It is easy to check that for every skew-Hamiltonian matrix $W\in\sHam$ there is a Hamiltonian matrix $H\in\Ham$
(and for every $H\in\Ham$ there is $W\in\sHam$) such that
$$W=\imath H.$$
Hence, the results given in the previous section can be applied here in a straightforward way.

%%%%%%%%%%%%%%%%%%%%%%%%%%%%%%%%%%%%%%%%%%%%%%%
%%%%%%%%%%%%%%%%%%%%%%%%%%%%%%%%%%%%%%%%%%%%%%%
%%%%%%%%%%%%%%%%%%%%%%%%%%%%%%%%%%%%%%%%%%%%%%%
\section{Structured canonical form for normal per(skew)-Hermitian  matrices under unitary perplectic similarity transformation}\label{sec:structured_perp}
The perskew-Hermitian matrices form a Lie algebra, while the per-Hermitian matrices form a Jordan algebra associated with $x^HFy.$ The perplectic matrices form the automorphism group associated with the skew-Hermitian sesquilinear form $x^HFy.$
Unlike in the (skew-)Hamiltonian case we did not find a suitable canonical form in the literature from which
a canonical form for normal and per(skew)-Hermitian matrices can be deduced. Thus we will directly
state and proof such a form.
Structure-preserving Jacobi-type algorithms for computing the structured canonical forms derived here  can be derived
adapting the ideas from Section \ref{sec:Jac_Ham}. We refrain from giving details.
%%%%%%%%%%%%%%%%%%%%%%%%%%%%%%%%%%%%%%%%%%%%%%%
\subsection{Canonical form for normal per-Hermitian matrices }
A matrix $M\in\mathbb{C}^{2n\times2n}$ is  per-Hermitian if
$(FM)^H=FM,$ or equivalently, $M^H=FMF$. If we write $M$ as a $2 \times 2$ block matrix of $n \times n$ blocks, it is easy to see that
\[
M = \left[ \begin{array}{cc}
M_{11} & M_{12}\\
M_{21} & FM_{11}^HF
\end{array}\right], \qquad\textnormal{where } (FM_{12})^H = FM_{12},  \quad (FM_{21})^H = FM_{21}.
\]
The elements of the antidiagonal of $M_{12}$ and $M_{21}$ have to be real.

As $F = F^H = F^{-1},$ we have $\sigma(M^H) = \sigma(FMF) = \sigma(M) = \sigma(M^T) =\overline{\sigma(M)}.$
Thus, eigenvalues with nonzero imaginary part appear in pairs $(\lambda,\overline{\lambda}).$ There is no restriction on the algebraic
multiplicity of real eigenvalues, in particular, it can be odd. Just the number of all real eigenvalues counted with multiplicity has to be even.

Obviously, a diagonal per-Hermitian matrix $\hat{M}$ has to be of the form
\[
\hat{M} = \left[ \begin{array}{cc}
\Lambda & 0\\
0 & F\Lambda^HF
\end{array}\right]
\]
with a diagonal matrix $\Lambda = \diag(\lambda_1, \ldots, \lambda_n),$ and, thus, $F\Lambda^HF =
\diag(\overline{\lambda}_n, \ldots, \overline{\lambda}_1).$

Hence, for any normal per-Hermitian matrix $M$ there exists a unitary matrix $U$ which diagonalizes $M$ such that
\[ D=U^HMU = \diag(D(\lambda_1) , D(\overline{\lambda_1}), \ldots, D(\lambda_t), D(\overline{\lambda_t}), D(\mu_1), \ldots, D(\mu_s))
\]
where
\begin{eqnarray*}
D(\lambda_j) &=& \left[\begin{array}{ccc}
\lambda_j & \\
& \ddots \\
&& \lambda_j
\end{array}\right] \in \mathbb{C}^{m_{\lambda_j} \times m_{\lambda_j}}, \quad
\lambda_j \in \mathbb{C}, j = 1, \ldots, t, \textnormal{ with Im}(\lambda_j) \neq 0,\\
D(\mu_j) &=& \left[\begin{array}{ccc}
\mu_j & \\
& \ddots \\
&& \mu_j
\end{array}\right]\in \mathbb{R}^{m_{\mu_j} \times m_{\mu_j}}, \quad \mu_j \in \mathbb{R}, j = 1, \ldots, s,
\end{eqnarray*}
and  $2n = 2\sum_{j=1}^tm_{\lambda_j} + \sum_{j=1}^s m_{\mu_j}$ with $\lambda_j \neq \lambda_i$ and $\mu_j \neq \mu_i,$ for $j \neq i.$
Moreover, $D(\overline{\lambda_j})  \in \mathbb{C}^{m_{\lambda_j} \times m_{\lambda_j}}.$
We will transform this unitary eigendecomposition into a similar decomposition of $M = V\widehat{D}V^H$
with a unitary perplectic $V$  and a matrix $\widehat{D}$ from which the eigenvalues can be read off immediately.

 Certainly, the columns of $U = [u_1~u_2~\cdots~u_{2n}]$ are orthonormal eigenvectors of $M.$
As $F$ is unitary we have that $D=U^HMU = (FU)^HM^H(FU),$ and $(FU)^HM(FU) = D^H.$
Hence, $FU$ is another unitary matrix which diagonalizes $M.$ Therefore, $u_j^HFu_k$ can only be nonzero in one of the following
situations:
\begin{enumerate}
\item $u_j$ is an eigenvector of $M$ for some $\lambda_\ell \in \mathbb{C}$ and $u_k$ is an eigenvector of $M$ for $\overline{\lambda_\ell}$ (or vice versa).
\item $u_j$ and $u_k$ are eigenvectors of $M$ for some $\mu_\ell \in \mathbb{R}.$
\end{enumerate}
In all other cases, $u_j^HFu_k=0.$ This implies that $U^HFU$ has the following form
$$U^HFU = \diag(
\left[\begin{smallmatrix} 0 & S(\lambda_1) \\ S(\lambda_1)^H & 0 \end{smallmatrix}\right],\ldots,
\left[\begin{smallmatrix} 0 & S(\lambda_t) \\ S(\lambda_t)^H & 0 \end{smallmatrix}\right],
 S(\mu_1), \ldots, S(\mu_s)).$$
As $U^HFU$ is unitary, all blocks $S(\lambda_j)$ and $S(\mu_j)$ have to be unitary. Moreover, as $U^HFU$ is Hermitian, the blocks
$S(\mu_j)$ are Hermitian.

Let $T \in \mathbb{C}^{2n \times 2n}$ denote the matrix
\[
T= \diag(S(\lambda_1), I_{m_{\lambda_1}}, \ldots, S(\lambda_t), I_{m_{\lambda_t}}, I_{2r})
\]
with $2r=\sum_{j=1}^s m_{\mu_j}.$ Clearly, $T$ is unitary. With $U_1 = UT$ we obtain
\[
U_1^HFU_1 = \diag( \left[\begin{smallmatrix}0 & I_{m_{\lambda_1}} \\ I_{m_{\lambda_1}} & 0 \end{smallmatrix}\right], \ldots,
\left[\begin{smallmatrix}0 & I_{m_{\lambda_t}} \\ I_{m_{\lambda_t}} & 0 \end{smallmatrix}\right], S(\mu_1), \ldots, S(\mu_s)),
\]
and
\[ U_1^HMU_1= D.\]
Next, we consider the blocks $S(\mu_j), j = 1, \ldots, s.$ Each such block is unitary and Hermitian.
Thus, all these blocks have eigenvalues $\pm 1$ (as $U_1^HFU_1$ is similar to $F$). Therefore, for each block $S(\mu_j)$ there exists a unitary matrix $V_j \in \mathbb{C}^{m_{\mu_j} \times m_{\mu_j}}$ such that $V_j^HS(\mu_j)V_j = \diag(\pm 1, \ldots, \pm 1)= \hat{I}_{\mu_j}$ for some combination of $+1$ and $-1$ entries.
Let
\[
V = \diag( I_{2c}, V_1, \ldots,  V_s) \in \mathbb{C}^{2n \times 2n}
\]
where $c = \sum_{j=1}^t m_{\lambda_j}.$ By construction, $V$ is unitary.
With $U_2= U_1V = UTV$ we obtain
\[
U_2^HFU_2 = \diag(\left[\begin{smallmatrix}0 & I_{m_{\lambda_1}} \\ I_{m_{\lambda_1}} & 0 \end{smallmatrix}\right], \ldots,
\left[\begin{smallmatrix}0 & I_{m_{\lambda_t}} \\ I_{m_{\lambda_t}} & 0 \end{smallmatrix}\right],  \hat{I}_{m_{\mu_1}}, \ldots,
\hat{I}_{m_{\mu_s}}),
\]
and
\[ U_2^HMU_2= D.\]
As $F$ has exactly $n$ eigenvalues $+1$ and $n$ eigenvalues $-1$ and the blocks $\left[\begin{smallmatrix} 0 & I_{m}\\
I_{m} & 0\end{smallmatrix}\right]$ have exactly $m$ eigenvalues $+1$ and $m$ eigenvalues $-1$, the diagonal block
\[
J = \diag(\hat{I}_{m_{\mu_1}}, \ldots, \hat{I}_{m_{\mu_s}}) \in \mathbb{R}^{2r \times 2r}
\]
 has to have exactly $r$ eigenvalues/diagonal entries  $+1$ and $r$ eigenvalues/diagonal entries $-1$
(recall that $2r=\sum_{j=1}^s m_{\mu_j}$).
These $\pm 1$ diagonal entries can be reordered by a simple permutation $\widetilde{P} \in \mathbb{R}^{2r \times 2r}$
such that
\[
\widetilde{P}^TJ\widetilde{P} =  \diag(+1, -1, +1, -1, \ldots, +1, -1).
\]
Let
\[
P_1 = \diag(I_{2c}, \widetilde{P})\in \mathbb{R}^{2n \times 2n}.
\]
As $\widetilde{P}$ and $P_1$ are unitary, the matrix $U_3 = U_2 P_1 = UTVP_1$ is unitary. Moreover,
\[
U_3^HFU_3 = \diag(\left[\begin{smallmatrix}0 & I_{m_{\lambda_1}} \\ I_{m_{\lambda_1}} & 0 \end{smallmatrix}\right], \ldots,
\left[\begin{smallmatrix}0 & I_{m_{\lambda_t}} \\ I_{m_{\lambda_t}} & 0 \end{smallmatrix}\right],  +1, -1, +1, -1, \ldots,
+1, -1),
\]
and
\[ U_3^HMU_3= D_1,\]
where $D_1$ is still diagonal, but not necessarily the same as $D,$ some of the real eigenvalues may have swapped places.
Next, we will make use of the fact that $\left[\begin{smallmatrix} +1 & 0\\ 0 & -1\end{smallmatrix}\right]$  and
$F_2 = \left[\begin{smallmatrix} 0 & 1\\ 1 & 0\end{smallmatrix}\right]$ are unitarily similar. Thus, there exists a unitary matrix
$Z \in \mathbb{C}^{2 \times 2}$ such that $Z^H\left[\begin{smallmatrix} +1 & 0\\ 0 & -1\end{smallmatrix}\right]Z = F_2.$
Let $\widetilde{W} = \diag(Z, \ldots, Z) \in \mathbb{C}^{2r \times 2r}$ and
\[
W = \diag(I_{2c}, \widetilde{W})\in \mathbb{C}^{2n \times 2n}.
\]
Now, as $\widetilde{W}$ and $W$ are unitary, the matrix $U_4 = U_3 W = UTVP_1W$ is unitary. Moreover,
\[
U_4^HFU_4 = \diag(\left[\begin{smallmatrix}0 & I_{m_{\lambda_1}} \\ I_{m_{\lambda_1}} & 0 \end{smallmatrix}\right], \ldots,
\left[\begin{smallmatrix}0 & I_{m_{\lambda_t}} \\ I_{m_{\lambda_t}} & 0 \end{smallmatrix}\right],  F_2,  F_2, \ldots,
F_2),
\]
and
\[ U_4^HMU_4= D_2 = \diag( D(\lambda_1), D(\overline{\lambda_1}), \ldots, D(\lambda_t), D(\overline{\lambda_t}),
X_1, X_2, \ldots, X_r)
\]
for some $X_j = \left[\begin{smallmatrix}a_j & b_j\\ b_j & a_j\end{smallmatrix}\right] \in \mathbb{R}^{2 \times 2}, j = 1, \ldots, r.$
Next, observe that $U_4^HFU_4$ is permutationaly similar to $F.$
This can be seen as follows. First, there is a permutation matrix $\widetilde{P}_2 \in \mathbb{R}^{2r \times 2r}$ such that $X = \widetilde{P}_2^T\diag(X_1, X_2, \ldots, X_r) \widetilde{P}_2$
is real, symmetric and per-symmetric
\begin{equation}\label{eq:X}
X =
\left[\begin{array}{ccccccccccccc}
a_1 &&&&&&&&&&b_1\\
& a_2 &&&&&&&&b_2&\\
&&\ddots &&&&&&\iddots\\
&&&& a_{s-1} &&&b_{s-1}\\
&&&&& a_s& b_s &\\
&&&&& b_s & a_s & \\
&&&&b_{s-1}&&&a_{s-1}\\
&&\iddots&&&&&&\ddots\\
&b_2&&&&&&&&a_2\\
b_1 &&&&&&&&&&a_1
\end{array}\right]
\end{equation}
and
\[
\widetilde{P}_2^T\diag(F_2, F_2, \ldots, F_2) \widetilde{P}_2 = F_{2r}.
\]
Moreover, there is a permutation matrix
$\widehat{P}_2\in \mathbb{R}^{2c \times 2c}$ which reorders the diagonal blocks $D(\lambda_j)$ and
$D(\overline{\lambda_j})$ in $U_4^HMU_4$ such that $P_2 =\diag(\widehat{P}_2,\widetilde{P}_2)$ achieves
\[
P_2^TU_4^HMU_4P_2 = \diag (D(\lambda_1), D(\lambda_2), \ldots, D(\lambda_t), D(\overline{\lambda_1}), D(\overline{\lambda_2}), \ldots, D(\overline{\lambda_t}), X)
\]
and
\[
P_2^TU_4^HFU_4P_2 = \diag(F_{2c},  F_{2r}).
\]
Finally, there is a $2n \times 2n$ permutation matrix $P_3$
such that
\[ P_3^TP_2^TU_4^HFU_4P_2P_3 = F_{2n}.\]
This permutation reorders the diagonal blocks of $P_2^TU_4^HMU_4P_2$ as follows
\[
P_3^TP_2^TU_4^HMU_4P_2P_3 = \diag (D(\lambda_1), D(\lambda_2), \ldots, D(\lambda_t), X,D(\overline{\lambda_1}), D(\overline{\lambda_2}), \ldots, D(\overline{\lambda_t})).
\]
Our findings can be summarized as follows (see also the discussion in \cite[Section 10.2]{Sal19})
\begin{Theorem}\label{theo:perHNormal}
Let $M \in \mathbb{C}^{2n \times 2n}$ be normal and per-Hermitian;  $M \in \N_{2n} \cap \perH.$
Then there exists a unitary and perplectic matrix $U\in \mathbb{C}^{2n \times 2n}$ such that
\[
U^HMU = \left[\begin{array}{ccc}
D\\ &X\\ &&FD^HF
\end{array}\right] \in \perH,
\]
where $D$ is a diagonal matrix whose diagonal entries have nonzero imaginary parts and $X$ is a real-valued matrix of the form \eqref{eq:X}. Moreover, $X$ is symmetric as well as per-symmetric. In particular, all eigenvalues of $X$ are real.
\end{Theorem}

%%%%%%%%%%%%%%%%%%%%%%%%%%%%%%%%%%%%%%%%%%%%%%%
\subsection{Canonical form for normal  perskew-Hermitian matrices under unitary perplectic similarity transformation}
A matrix $K\in\mathbb{C}^{2n\times2n}$ is  perskew-Hermitian if
$(FK)^H=-FK,$ or equivalently, $K^H=-FKF$. 
It is easy to check that for every perskew-Hermitian matrix $K\in\sperH$ there is a per-Hermitian matrix $M\in\perH$
(and for every $M\in\perH$ there is $K\in\sperH$) such that
$$K=\imath M.$$
Hence, the results given in the previous section can be applied here in a straightforward way.

%%%%%%%%%%%%%%%%%%%%%%%%%%%%%%%%%%%%%%%%%%%%%%%
%%%%%%%%%%%%%%%%%%%%%%%%%%%%%%%%%%%%%%%%%%%%%%%
%%%%%%%%%%%%%%%%%%%%%%%%%%%%%%%%%%%%%%%%%%%%%%%
\section{Algorithm for computing the canonical form of normal (skew-)Hamiltonian matrices}\label{sec:Jac_Ham}
The algorithm for diagonalizing normal matrices as suggested in \cite{GolH59} is based on 
 the idea  of writing the normal matrix $A$ under consideration as a sum of its Hermitian
part $B= \frac{1}{2}(A + A^H)$ and skew-Hermitian part $C= \frac{1}{2}(A - A^H).$ The unitary matrix $U$ which diagonalizes $A$ also diagonalizes $B$ and $C.$
If the order of the eigenvalues on the diagonal of $D= U^HAU$  is fixed, $U$ is essentially unique  in the following sense.
For single eigenvalues of $A,$ the corresponding eigenvector in $U$ is unique up to 
multiplication with a complex sign $e^{\imath \alpha}, \alpha \in \mathbb{R}.$ For an eigenvalue 
of $A$ with multiplicity $r_j$, the corresponding $r_j$  eigenvectors in $U$ form an invariant subspace of $A.$ 
Now, assume that a unitary matrix $V \in \mathbb{C}^{n \times n}$ diagonalizing $B$ 
has been found by, e.g., the classical Jacobi method \cite{Jac1846,GolVL12} adapted for Hermitian matrices \cite{Hen58}.
Then, $V^HCV$ is not necessarily diagonal. But the only nonzero off-diagonal
elements in $V^HCV$ are in those positions determined by the equal eigenvalues of $B.$ 
The  simple Givens rotations can be employed to complete the diagonalization of $V^HCV.$ They commute with 
the diagonal matrix $V^HBV$ and therefore leave it invariant. 

This idea can be adapted for all four types of structured normal matrices considered here. A main ingredient for such
an algorithm is the structure-preserving diagonalization of the Hermitian part of the matrix under consideration.
As such a Jacobi-type algorithm is only known to us in the Hamiltonian case, we will sketch how to adapt the idea above 
for the Hamiltonian case and present a proof of concept. 
In the literature one can find two Jacobi-type algorithms for computing the Hamiltonian Schur form of Hamiltonian matrices
 with no purely imaginary eigenvalues \cite{Bye90,BunF97}. These may be used to diagonalize Hermitian Hamiltonian
matrices. Please note that these algorithms can not be employed here directly, as they can not be used  for diagonalizing general normal Hamiltonian matrices which may have purely imaginary eigenvalues. 

Let $H\in\N_{2n} \cap \Ham.$ 
 Then its Hermitian part
$B = \frac{1}{2}(H+H^H)$ is Hermitian and Hamiltonian, while its skew-Hermitian part 
$C = \frac{1}{2}(H-H^H)$ is skew-Hermitian and Hamiltonian. A unitary and symplectic matrix $Z$ which
transforms $H$ to its Hamiltonian Schur form \eqref{eq:HamCanonical1} yields
\begin{align}
Z^HBZ &= \begin{bmatrix}
\re(D_1) & 0 & 0 &0\\
0 & 0 & 0 & 0\\
0 & 0 & -\re(D_1) & 0 \\
0 & 0& 0 & 0
\end{bmatrix},   \label{eq:B}\\
Z^HCZ &=  \begin{bmatrix}
\imath\im(D_1) & 0 & 0 &0\\
0 & D_2 & 0 & D_3\\
0 & 0 & \imath\im(D_1) & 0 \\
0 & -D_3 & 0 & D_2
\end{bmatrix} \label{eq:C}
\end{align}
where $D_1, D_2, D_3$ are as in \eqref{eq:HamCanonical1}.
$Z$ is essentially unique. For single eigenvalues of $B,$ the corresponding eigenvector in $Z$ is unique up to 
multiplication with a complex sign $e^{\imath \alpha}, \alpha \in \mathbb{R}.$ For an eigenvalue of $B$ with multiplicity $r_j$, 
the corresponding $r_j$  eigenvectors in $Z$ form an invariant subspace of $B.$ The order of these eigenvalues
on the diagonal of $Z^HBZ$ is arbitrary.
 Assume we have constructed a unitary and symplectic matrix $S$ such that 
\begin{equation}\label{eq:B2}
S^HBS = \begin{bmatrix}
\Lambda & 0 & 0 &0\\
0 & 0 & 0 & 0\\
0 & 0 & -\Lambda & 0 \\
0 & 0& 0 & 0
\end{bmatrix}
\end{equation}
with $\Lambda = \diag(\lambda_1, \ldots, \lambda_{n_1})\in \mathbb{R}^{n_1 \times n_1}.$  
Assume for ease of notion that the nonzero eigenvalues of $B$ in \eqref{eq:B} and \eqref{eq:B2} are ordered such that
\[
\lambda_1 = \re((D_1)_{11}) \geq  \lambda_2 = \re((D_1)_{22}) \geq \cdots \geq \lambda_{n_1} = \re((D_1)_{{n_1}{n_1}}),
\]
that is, $Z^HBZ = S^HBS.$
For all single eigenvalues of $B,$ the corresponding eigenvectors in $S$ and $Z$ are identical up to 
multiplication with a complex sign. For any multiple eigenvalue of $B$ the corresponding eigenvectors in $S$ and $Z$ 
form the same invariant subspace. Thus, in case $B$ has $p$ distinct nonzero eigenvalues $\mu_j$ with multiplicities $r_j$
we have
\[
S^HBS = Z^HBZ = \diag(\Lambda_1, \ldots, \Lambda_p,0,-\Lambda_1, \ldots, -\Lambda_p,0), \quad \Lambda_j = \mu_j I_{r_j} \in \mathbb{C}^{r_j \times r_j}
\]
and
\[
S = Z \left[W \oplus W\right], \qquad W= \left[ W_1 \oplus W_2 \oplus \cdots \oplus W_p \oplus W_{n-n_1}\right]
\]
for unitary $W_j \in \mathbb{C}^{r_j \times r_j}, j = 1, \ldots, p$ and $\sum_{\ell=1}^p r_\ell = n_1.$

From this, we see that $S^HCS$ will not be as in \eqref{eq:C} as 
\begin{align*}
S^HCS &= \left[W \oplus W\right]^HZ^HCZ\left[W \oplus W\right]\\
&=\left[W \oplus W\right]^H
\begin{bmatrix}
\imath\im(D_1) & 0 & 0 &0\\
0 & D_2 & 0 & D_3\\
0 & 0 & \imath\im(D_1) & 0 \\
0 & -D_3 & 0 & D_2
\end{bmatrix} 
\left[W \oplus W\right]\\
&= 
\begin{bmatrix}
C_1 & 0 &0 &0\\
0 & C_2 & 0 & C_3\\
0 & 0 & C_1 & 0 \\
0 & -C_3 & 0 & C_2
\end{bmatrix} 
\end{align*}
with $C_1= -C_1^H, C_2= -C_2^H, C_3 = C_3^H.$
The blocks $C_2$ and $C_3$ may
have nonzero entries at every position. Moreover,
if the $(k,k)$ and the
$(j,j)$ entries in $\re(D_1)$ are equal, then in the block $C_1$ 
there may be entries at the positions $(k,j)$ and $(j,k)$. 

The submatrix $\left[\begin{smallmatrix}C_2 & C_3\\ -C_3 & C_2\end{smallmatrix}\right]$ is a
normal skew-Hermitian Hamiltonian matrix with purely imaginary eigenvalues. It needs to be diagonalized
by a unitary symplectic matrix. This can be done independently of what needs to be done to 
diagonalize the submatrix $\left[\begin{smallmatrix}C_1 & 0\\ 0 & C_1\end{smallmatrix}\right].$ 

As we do not know of a structure-preserving Jacobi-type algorithm for  skew-Hermitian Hamiltonian matrices, we suggest to transform the
 submatrix $\left[\begin{smallmatrix}C_2 & C_3\\ -C_3 & C_2\end{smallmatrix}\right]\in \mathbb{C}^{2(n-n_1) \times 2(n-n_1)}$  
into the desired form as described in
the proof of Theorem \ref{theo:HamSchurNormal}. That is, first unitary matrices $V_1$ and $V_2$
are determined which diagonalize the skew-Hermitian matrices $C_2+\imath C_3$ and
$C_2-\imath C_3 \in \mathbb{C}^{(n-n_1) \times (n-n_1)}.$  Then the unitary symplectic matrix
$\hat{S} = QVQ^H$ with $Q = \frac{1}{\sqrt{2}} \left[ \begin{smallmatrix} I & \imath I\\ \imath I & I\end{smallmatrix}\right]$
and $V = \diag(V_1, V_2)$ diagonalizes each block of $\left[\begin{smallmatrix}C_2 & C_3\\ -C_3 & C_2\end{smallmatrix}\right];$
$\hat{S} = \left[\begin{smallmatrix}\hat{S}_1 & \hat{S_2} \\ -\hat{S}_2 & \hat{S}_1\end{smallmatrix}\right] 
= \left[\begin{smallmatrix}V_1+V_2 & \imath(V_2-V_1) \\ -\imath(V_2-V_1) & V_1+V_2\end{smallmatrix}\right].$ 
As any skew-Hermitian matrix $G$ can be expressed as $G = \imath F$ for a Hermitian matrix $F,$
the unitary matrices $V_1$ and $V_2$ can be computed using the classical Jacobi method for Hermitian matrices.

Finally, let us consider $\left[\begin{smallmatrix}C_1 & 0\\ 0& C_1\end{smallmatrix}\right]
\in \mathbb{C}^{2n_1 \times 2n_1}.$  Assume that the entries $(k,j)$ and $(j,k)$ in $C_1$ are nonzero (and, hence,
the diagonal entries $(k,k)$ and $(j,j)$ of $\re(D_1)$ are equal). Our goal is to annihilate the entries $(C_1)_{kj}$ and 
$-(C_1)_{kj}$ without altering $\Lambda.$
Thus,   consider the $2 \times 2$ subproblem of $C_1,$
\[
\widehat{C} = \begin{bmatrix} 
(C_1)_{kk} &(C_1)_{kj} \\ -(C_1)_{kj}  & (C_1)_{j,j} 
\end{bmatrix},
\]
as well as
 the corresponding $2 \times 2$ subproblem of $\Lambda,$
\[
\widehat{B}=\begin{bmatrix} 
 \re(D_1)_{kk} &0\\ 0  &\re(D_1)_{jj} 
\end{bmatrix},
\]
with $\re(D_1)_{kk} =\re(D_1)_{jj}.$ As $\widehat{C}$ is normal, there exists a unitary matrix 
\[\widehat{G} =
\begin{bmatrix} \cos x & -e^{-\imath \alpha} \sin x\\ e^{\imath \alpha}\sin x & \cos x \end{bmatrix} =
\begin{bmatrix}  c& -s \\ \bar{s} & c\end{bmatrix}
\]
with $x, \alpha \in \mathbb{R}$ such that
$\widehat{G}^H\widehat{C}\widehat{G}$ is diagonal. Moreover,
\[
 \widehat{G}^H\widehat{B}\widehat{G} = \diag( \re(D_1)_{kk}, \re(D_1)_{jj}) = \widehat{B}
\]
holds. Now, let $G_{kj}(c,s) \in \mathbb{C}^{n \times n}$  be a unitary Givens rotation
\begin{equation}\label{eq:givens}
G_{kj}(c,s) = \left[
\begin{array}{ccccc}
I_{k-1}\\
& c && -s\\
&&I_{j-k-1}\\
& \bar{s} && c\\
&&&& I_{n-j}
\end{array}\right] \in \mathbb{C}^{n \times n},
\end{equation}
where $k$ and $j$ describe the positions of $c$ and $s,$ $1 \leq k< j \leq n.$ 
Further, let $R_{kj}(c,s)$ be the direct sum $G_{kj} (c,s)\oplus G_{kj}(c,s) \in \mathbb{C}^{2n \times 2n}$ with $ G_{kj}(c,s) \in \mathbb{C}^{n \times n}.$ This yields a unitary and symplectic transformation matrix, called a
symplectic direct sum embedding.
Then $R_{kj}$ 
 annihilates the entries $(k,j)$ and $(j,k)$ in $C_1.$ This will not alter $B.$
With the help of  these transformations, $C_1$ can be diagonalized while $B$ remains diagonal. 
Denote the unitary matrix which diagonalizes $C_1$ by $\hat{S}_3 \in \mathbb{C}^{n_1 \times n_1}.$ Then $\diag(\hat{S}_3, \hat{S}_3)$ is an unitary symplectic matrix which diagonalizes $\diag(C_1,C_1)$
and does not change $\diag(\Lambda,\Lambda).$
 In summary, with the unitary symplectic matrix
\[
ZT= \begin{bmatrix}
\hat{S}_3\\ & I_{n-n_1}&\\ &&\hat{S}_3\\&&&I_{n-n_1}
\end{bmatrix}\begin{bmatrix}
I_{n_1}\\ & \hat{S}_1&&\hat{S}_2\\ &&I_{n_1}\\&-\hat{S}_2&&\hat{S}_1
\end{bmatrix}
\]
we have
\[
T^HS^HAST = \begin{bmatrix}
D_1 \\ & D_2 &&D_3\\
&&-D_1^H\\
&-D_3&&D_2
\end{bmatrix}.
\]

In the very first step of the approach described here, a unitary and symplectic matrix $S$ needs to be found 
which diagonalizes $B$ as in \eqref{eq:B2}. For this step,
we can make use of one of the two Jacobi-type algorithms for computing the Hamiltonian Schur form of Hamiltonian matrices
 with no purely imaginary eigenvalues \cite{Bye90,BunF97}. Hermitian Hamiltonian matrices have just real eigenvalues. Thus, the algorithms can be applied. We will briefly review these algorithms.

The Hamiltonian-Jacobi algorithm proposed in \cite{Bye90}  is based on unitary symplectic Givens transformations
$G_{j,n+j}(c,s)$ as in \eqref{eq:givens}.
 It is noted in \cite{Bye90}
that the algorithm converges too slowly (for use on conventional serial computers) as its convergence properties
are like the ones discussed in \cite{Ste85}.  Thus this algorithm will not be considered here any further.

The basic idea of the algorithm presented in  \cite{BunF97} is to consider
$4 \times 4$ (Hermitian) Hamiltonian submatrices $H_{ij}$ of $H = \left[\begin{smallmatrix}A & F\\ G & -A^H\end{smallmatrix}\right],$
\[
H_{ij} = \left[\begin{array}{cc|cc}
a_{ii} & a_{ij} & f_{ii} & f_{ij}\\
{a_{ji}} & a_{jj} & \bar{f}_{ij} & f_{jj}\\ \hline
g_{ii} & g_{ij} & -\bar{a}_{ii} & -\bar{a}_{ji}\\
\bar{g}_{ij} & g_{jj} & -\bar{a}_{ij} & -\bar{a}_{jj}
\end{array} \right].
\]
These submatrices are transformed as far as possible to their (Hermitian) Hamiltonian Schur form by a symplectic unitary transformation. This transformation can be constructed directly in a finite number of steps once one eigenvalue of $H_{ij}$ is known. 
In case, $H_{ij}$ has only purely imaginary eigenvalues, it is left unchanged as the transformation matrix is chosen to be
the identity. 
As discussed in \cite{BunF97}, for Hermitian Hamiltonian matrices $H = \left[ \begin{smallmatrix} A & G\\ G &-A\end{smallmatrix}\right],$ $A=A^H,$ $A=G^H$
this approach for computing the Hamiltonian Schur form is equivalent to Kogbetliantz's Jacobi-type
algorithm for computing the SVD of $A+\imath G.$
Quadratic convergence of Kogbetliantz's Jacobi-type algorithm has been proven in \cite{PaiVD86} in case that $A+\imath G$ does not have close singular values.

This concludes our sketch of the algorithm to compute the structured canonical form of a normal Hamiltonian matrix.

\subsection{Convergence of the proposed algorithm}
The algorithm proposed above is convergent. In order to see this,
let us consider the first step of our computation on the Hermitian part. Under mild assumptions this step is quadratically convergent. Denote the accumulated
transformation matrix used in the $j$th sweep of the algorithm by $S_j.$ Then, there is a scalar $k_1\in \mathbb{N}$  such that after performing $k_1$ sweeps
we have for $\widetilde{S} = S_1S_2\cdots S_{k_1}$
\begin{align*}
A^{(k_1)} &= \widetilde{S}^HA\widetilde{S} =   \widetilde{S}^HB\widetilde{S} + \widetilde{S}^HC\widetilde{S} \\
&=
\left( \begin{bmatrix}
\re(D_1) & 0 & 0 &0\\
0 & 0 & 0 & 0\\
0 & 0 & -\re(D_1) & 0 \\
0 & 0& 0 & 0
\end{bmatrix} + E\right) +
\left( \begin{bmatrix}
C_1 & 0 & 0 &0\\
0 & C_2 & 0 & C_3\\
0 & 0 & C_1 & 0 \\
0 & -C_3& 0 & C_2
\end{bmatrix} + F\right)
\end{align*}
with $\lim_{k_1\rightarrow \infty}\|E\|_F=0$ and $\lim_{k_1\rightarrow \infty}\|F\|_F=0.$
The second step of our computation on the skew-Hermitian part uses the Jacobi algorithm for Hermitian matrices. We can view this as performing sweeps on $\widetilde{S}^HC\widetilde{S}.$ Denote the transformation matrix used in the $j$th sweep of the algorithm $T_j.$ Then,
there is a scalar $k_2\in \mathbb{N}$  such that after performing $k_2$ sweeps
we have for $\widetilde{T} = T_1T_2\cdots T_{k_2}$
\begin{align*}
&\widetilde{T}^H\widetilde{S}^HA\widetilde{S}\widetilde{T} =   
\widetilde{T}^H\widetilde{S}^HB\widetilde{S}\widetilde{T}  + \widetilde{T}^H\widetilde{S}^HC\widetilde{S}\widetilde{T}  \\
&\quad =
\left( \begin{bmatrix}
\re(D_1) & 0 & 0 &0\\
0 & 0 & 0 & 0\\
0 & 0 & -\re(D_1) & 0 \\
0 & 0& 0 & 0
\end{bmatrix} + \widetilde{T}^HE\widetilde{T}\right) \\
&\qquad\qquad+
\left( \begin{bmatrix}
\imath \im(D_1) & 0 & 0 &0\\
0 & D_2 & 0 & D_3\\
0 & 0 & \imath \im(D_1) & 0 \\
0 & -D_3& 0 & D_2
\end{bmatrix} + \widetilde{T}^HF\widetilde{T}\right) + G
\end{align*}
with $\lim_{k_2\rightarrow \infty}\|\widetilde{T}^HE\widetilde{T}\|_F=0,$ $\lim_{k_2\rightarrow \infty}\|\widetilde{T}^HF\widetilde{T}\|_F=0$ and $\lim_{k_2\rightarrow \infty}\|G\|_F=0.$ In summary we have global convergence 
\[
\lim_{k_1, k_2\rightarrow \infty} \widetilde{T}^H\widetilde{S}^HA\widetilde{S}\widetilde{T} = 
 \begin{bmatrix}
D_1 & 0 & 0 &0\\
0 & D_2 & 0 & D_3\\
0 & 0 & D_1 & 0 \\
0 & -D_3& 0 & D_2
\end{bmatrix}.
\]

\subsection{Numerical Experiments}
The algorithm proposed here has been implemented in MATLAB R2019a. We have used
 the implementation of the Jacobi-type method for computing the Hamiltonian Schur form as given in \cite{BunF97}. 
Moreover, the algorithm from \cite{GolH59} for diagonalizing normal matrices has been implemented using the implementation of the Jacobi method for Hermitian matrices from \cite{web}.

The algorithm for computing the canonical form \eqref{eq:HamCanonical1} of normal Hamiltonian matrices   exhibits the same numerical behaviour as the one for diagonalizing normal matrices from \cite{GolH59}. We will demonstrate this by considering a normal Hamiltonian matrix $H$ for which its eigenvalues as well as the eigenvalues of its Hermitian
part $B = (H+H^H)/2$ are all distinct. In this case the unitary symplectic matrix $S$ in \eqref{eq:B2} will (at least theoretically)
diagonalize $C.$ Similarly, the (unstructured) unitary matrix $U$ which diagonalizes $B$ will (at least theoretically)
diagonalize $C.$ 
%As the following small example shows in practice this is not always the case.
 
For this example, we have constructed a (random) normal Hamiltonian matrix $H\in \mathbb{C}^{2n \times 2n}$ in Schur form \eqref{eq:HamCanonical1} by choosing $n_1 = n.$ In particular, we use the following code to set up $H$\\
\texttt{D1 = diag(randn(n,1))+1i*diag(randn(n,1));} \texttt{D = [D1 zeros(n,n); zeros(n,n) -D1'];}\\
\texttt{V1 = orth(randn(n,n)+1i*randn(n,n);} \texttt{V2 = orth(randn(n,n)+1i*randn(n,n);} \\
\texttt{W = [V1+V2  1i*(V2-V1);  -1i*(V2-V1)  V1+V2]/2;} \texttt{H=SW*D*W;}  \texttt{B = (H+H')/2;} \texttt{C = (H-H')/2;}

\begin{figure}[!htb]
\begin{center}
\includegraphics[width=.7\textwidth]{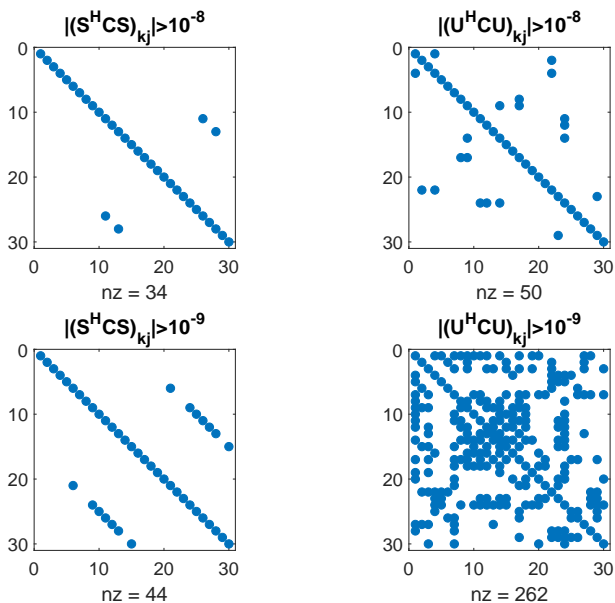}
\end{center}
\caption{Nonzero pattern of the skew-Hermitian part of a $30 \times 30$ normal Hamiltonian matrix $H$ which can be 
diagonalized by a unitary symplectic transformation,  structure-preserving algorithm on the left, unstructured algorithm on the right.}\label{fig1}
\end{figure}

%We use the algorithm from \cite{BunF97} to diagonalize the Hermitian part $B$ in a structure-preserving way.
%Let us denote the final unitary symplectic matrix which diagonalizes $B$ by $S.$
%We also use the algorithm from \cite{web} to diagonalize the Hermitian part $B.$ This algorithm is not 
%structure-preserving, it computes a unitary matrix $U$ which diagonalizes $B.$ 
Due to the set up of this example, $S^HCS$ and $U^HCU$ should be diagonal  for the unitary symplectic matrix $S$ which diagonalizes $B$ and the (unstructured) unitary matrix $U$ which diagonalizes $B.$
As can be seen in Fig. \ref{fig1}, in practice this is (unsurprisingly) only partially the case. Just working on $B$ and taking $B$ into account in determining when to stop the algorithms is not sufficient.
Both algorithms use different stopping criteria. The 
algorithm from \cite{BunF97}  (used to compute the unitary symplectic $S$ which diagonalizes $B$) stops when the relative error $\frac{\|B^{(j)}-\operatorname{diag}(B^{(j)})\|_F}{\|B\|_f}$ 
for the current iterate $B^{(j)}$ is less than
a given tolerance $tol.$ The algorithm from \cite{web} (used to compute the unitary $U$ which diagonalizes $B$) stops when the absolute values of all off-diagonal entries of the current iterate are less than $tol.$ 
The example shown in Fig. \ref{fig1} is of size $30 \times 30,$ $tol$ is chosen as $10^{-10}.$
The absolute error in the eigenvalues is of the order
of $10^{-14}$ for both algorithms. The structure-preserving algorithm stopped after six sweeps, while the
unstructured algorithm needed eight sweeps. 
In order to accommodate for the different stopping criteria, the spy plots in upper row of Fig. \ref{fig1}
show all nonzero elements which are larger than $100 \cdot tol,$
while the spy plots in the lower row of Fig. \ref{fig1} show all nonzero elements which are larger than $10 \cdot tol.$
 The displayed situation is typical. Using the structure-preserving 
algorithm, there are usually some entries on the diagonals of the $(1,2)$ and the $(2,1)$ block that remain, while using the unstructured
algorithm no particular pattern of the remaining nonzero elements can be observed. It is not possible to annihilate those unwanted entries by a simple sequence of symplectic direct sum embeddings in the case of the structure-preserving 
algorithm or of Givens rotations in the case of the unstructured algorithm. Those transformations would alter $S^HBS$ and $U^HBU.$ The only way to drive the entries in $S^HCS$ and $U^HCU$ to zeros is to use improve stopping criteria depending not just on $B$, but also on $C$, for the algorithms which determine $S$ and $U$ by just working on $B.$

\begin{figure}[!htb]
\begin{center}
\includegraphics[width=.7\textwidth]{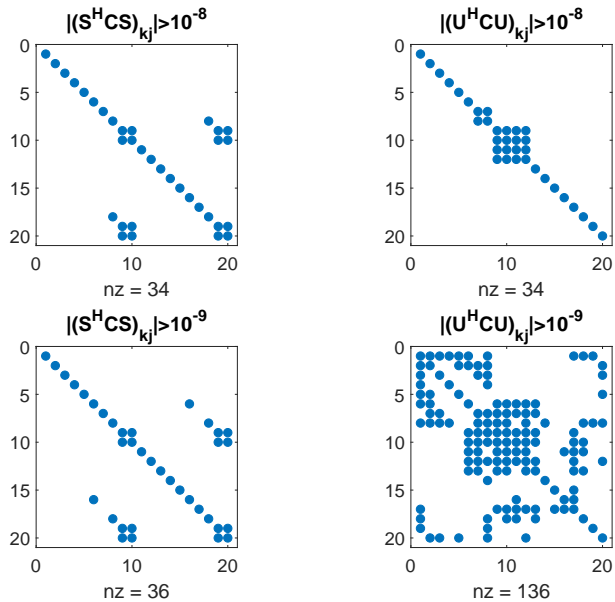}
\end{center}
\caption{Nonzero pattern of the skew-Hermitian part of a $30 \times 30$ normal Hamiltonian matrix $H$ with purely imaginary eigenvalues,  structure-preserving algorithm on the left, unstructured algorithm on the right.}\label{fig2}
\end{figure}

A similar observation can be made for normal Hamiltonian matrices with purely imaginary eigenvalues. The example used to
generate Fig. \ref{fig2} is a $20 \times 20$ normal Hamiltonian matrix with sixteen complex eigenvalues (with distinct nonzero real part) and four purely imaginary eigenvalues. As before, we use the algorithm from \cite{BunF97} to diagonalize the Hermitian part $B$ in a structure-preserving way. Denote the final unitary symplectic matrix which diagonalizes $B$ by $S.$ Then
$S^HCS$ will not be diagonal, the purely imaginary eigenvalues will results in $2 \times 2$ blocks in each of the four
$n \times n$ blocks of $S^HCS.$ This can be seen nicely in Fig. \ref{fig2}. But as before, there are additional unwanted
nonzero elements in $S^HCS.$ 
We also use the algorithm from \cite{web} to diagonalize the Hermitian part $B.$ Its effect on $C$ is displayed in the right
column in Fig. \ref{fig2}. As observed before, good stopping criteria for both algorithms must depend not just on $B$, but also on $C.$
A full analysis of this is beyond the scope of this paper.

\section{Concluding Remarks}
Normal and (skew-)Hamiltonian as well as normal and per(skew)-Hermitian matrices have been considered.
Structured canonical forms under unitary and suitable structure-preserving similarity transformations have been given.
A structure-preserving algorithm for computing the structured canonical form for normal Hamiltonian matrices 
has been presented. Some of its numerical properties have been discussed. It can easily be adapted for computing the canonical forms of the other normal structured matrices
considered here.

%%%%%%%%%%%%%%%%%%%%%%%%%%%%%%%%%%%%%%%%%%%%%%%

\section*{Acknowledgement}
The first author acknowledges financial support by DAAD Short-term grant and by Croatian Science Foundation under the project 3670.
The authors would like to thank Christian Mehl for discussions on Section \ref{sec:subsec1}.
%%%%%%%%%%%%%%%%%%%%%%%%%%%%%%%%%%%%%%%%%%%%%%%

\end{document}